\documentclass{amsart} 
\usepackage{amssymb,latexsym,amscd,hyperref}

\newtheorem{theorem}{Theorem}[section]

\newtheorem{e-proposition}[theorem]{Proposition}

\newtheorem{conjecture}[theorem]{Conjecture}
\newtheorem{e-definition}[theorem]{Definition\rm}

\def\og{\leavevmode\raise.3ex\hbox{$\scriptscriptstyle\langle\!\langle$~}}
\def\fg{\leavevmode\raise.3ex\hbox{~$\!\scriptscriptstyle\,\rangle\!\rangle$}}

\newcommand {\Q}{{\mathbb{Q}}}

\newcommand {\FF}{{\mathbb{F}}}
\newcommand{\Gal}      {\mathop{\rm {Gal}}}
\newcommand{\ind}        {{\mathop{\rm ind}}}

\newcommand{\eps}{\epsilon}
\newcommand{\Norm}{{\rm N}}

\begin{document}


\title[A counter example]{A counter example to Malle's conjecture on the asymptotics of discriminants}
\vspace{-2.6cm}


\author{J\"urgen Kl\"uners}
\email{klueners@mathematik.uni-kassel.de}
\address{Universit\"at Kassel, Fachbereich Mathematik/Informatik,
Heinrich-Plett-Str. 40, 34132 Kassel, Germany.}

\begin{abstract}
  In this note we give a counter example to a conjecture of Malle which
  predicts the asymptotic behaviour of the counting functions for field
  extensions with given Galois group and bounded discriminant. 
\end{abstract}
\vskip 0.5\baselineskip
\maketitle


\section{Introduction}
Let $G\le S_n$ be a finite transitive permutation group
and $k$ be a number field. We say that a finite extension
$K/k$ has Galois group $G$ if the normal closure $\hat K$ of $K/k$ has Galois
group isomorphic to $G$ and $K$ is the fixed field in $\hat K$ under a point
stabilizer of $G$. By abuse of notation we will write $\Gal(K/k)=G$ in this
situation. We let
$$Z(k,G;x):=|\left\{K/k : \Gal(K/k)=G,\ \Norm_{k/\Q}(d_{K/k})\le x\right\}|$$
be the number of field extensions of $k$ (inside a fixed algebraic closure
$\bar\Q$) of degree~$n$ with Galois group permutation isomorphic to $G$
(as explained above) and norm of the discriminant $d_{K/k}$ bounded above
by $x$. It is
well known that the number of extensions of $k$ with bounded norm of the
discriminant is finite, hence $Z(k,G;x)$ is finite for all $G$, $k$ and
$x\ge1$. 

Gunter Malle \cite{Ma4,Ma5} has given a precise conjecture about the asymptotic 
behaviour of the function $Z(k,G;x)$ for $x\rightarrow\infty$. In order
to state it, we need to introduce some group theoretic invariants of permutation
groups. 
\begin{e-definition}
  Let $1\ne G\leq S_n$ be a transitive subgroup acting on $\Omega=\{1,\ldots,n\}$.
  \begin{enumerate}
  \item For $g\in G$ we define $\ind(g):= n- \mbox{ the number of orbits of $g$ on }\Omega.$
  \item $\ind(G):=\min\{\ind(g): 1\ne g\in G\}.$
  \item $a(G):=\ind(G)^{-1}$.
  \end{enumerate}
\end{e-definition} 
Since all elements in a conjugacy class $C$ of $G$ have the same index
we can define $\ind(C)$ in a canonical way. The absolute Galois group
of $k$ acts on the set of conjugacy classes of $G$ via the action on the
columns of the complex character table of $G$. The orbits under this action
are called $k$-conjugacy classes.
\begin{e-definition}
  For a number field $k$ and a transitive subgroup $1\ne G\leq S_n$ we define:
  $$b(k,G):=|\{C : C\; k\mbox{-conjugacy class of minimal index }\ind(G)\}|.$$
\end{e-definition}
Now we can state the conjecture of Malle\cite{Ma5}, where $f(x) \sim g(x)$ is a 
notation for $\lim_{x\rightarrow\infty}f(x)/g(x) =1$.
\begin{conjecture}\label{con}(Malle)
  For all number fields $k$ and all transitive permutation groups $1\ne G$ there exists a 
  constant $c(k,G)>0$ such that
  $$Z(k,G;x) \sim c(k,G)x^{a(G)} \log(x)^{b(k,G)-1},$$
  where $a(G)$ and $b(k,G)$ are given as above.  
\end{conjecture}
  This conjecture is proved for abelian groups \cite{wr}. 
For all number fields $k$ and all nilpotent groups $G$ it is shown in \cite{KlMa2} that
$$\limsup_{x\rightarrow\infty} \frac{\log Z(k,G;x)}{\log x}\leq a(G).$$
If we furthermore assume that $G$ is in its regular representation, i.e.,
we count normal nilpotent number fields, we get:
$$\lim_{x\rightarrow\infty} \frac{\log Z(k,G;x)}{\log x}= a(G).$$
For more results see also the survey articles \cite{Bel,CoDiOl4}.

\section{The counter example}
\label{sec:counter}

We present a counter example for $k=\Q$ and the wreath product 
$G:=C_3\wr C_2 = C_3^2\rtimes C_2 \leq S_6$ of order 18 , where $C_n$ is
the cyclic group of order $n$.  In the following we count all field
towers $L/K/\Q$ such that $\Gal(L/K)=C_3$ and $\Gal(K/\Q)=C_2$.
Therefore the Galois group of $L/\Q$ is one of the groups $C_6,
S_3(6), G\leq S_6$, where $S_3(6)$ denotes the group $S_3$ in its
degree 6 representation.  Since $C_6$ is abelian we get from \cite{wr} that 
$$Z(\Q,C_6;x) \sim c(C_6) x^{1/3}.$$
From the Davenport-Heilbronn theorem \cite{DaHe} we know that 
$$Z(\Q,S_3(3);x) \sim c(S_3(3)) x.$$
Using the fact that the discriminant of the splitting field of an
$S_3$-extension is at least the square of the discriminant of the $S_3$-extension, we
easily get that $$Z(\Q,S_3(6);x)= O(x^{1/2}).$$
Since extensions with Galois group $S_3(6)$ are normal we can use a result in \cite[Prop. 2.8]{ElVe}
which states that $Z(\Q,S_3(6);x)= O_{\eps}(x^{3/8+\eps})$ for
all $\eps>0$. With a more careful analysis we are able to prove that
there are constants $c_1(S_3), c_2(S_3)>0$ such that
$$c_1(S_3)x^{1/3} \leq Z(\Q;S_3(6);x) \leq c_2(S_3)x^{1/3}\;\;
\mbox{ for }x \mbox{ large enough}.$$
We remark that the results for
$S_3(6)$ and $C_6$ are as conjectured since we get
$a(C_6)=a(S_3(6))=1/3$ and $b(\Q,C_6)=b(\Q,S_3(6))=1$.

Now we define the counting function corresponding to field towers $L/K/\Q$ as above:
$$\tilde Z(\Q,C_3\wr C_2;x) := |\left\{L/\Q\mid \exists K:\Gal(L/K)=C_3, [K:\Q]=
2, 
|d_L|\le x\right\}|.$$
We have two conjugacy classes of elements of order 3 in $G$ which have three fixed points. Considered
as $\Q$-conjugacy classes we get only one orbit. In number fields $k$ containing a primitive third root
of unity $\zeta_3$ we get two $k$-conjugacy classes of this type. Therefore $a(G)=1/2$ and $b(\Q,G)=1$.
Since the counting functions for $S_3(6)$ and $C_6$ have lower asymptotics we would expect that
$$\tilde Z(\Q,C_3\wr C_2;x) \sim Z(\Q,C_3\wr C_2;x) \sim c(G) x^{1/2}.$$
Certainly we get a lower estimate for $\tilde Z(\Q,C_3\wr C_2;x)$ if we only count the number fields
which contain a fixed quadratic subfield $K$. We choose $K=\Q(\zeta_3)$ and using 
$d_L=d_K^3 \Norm(d_{L/K})$ we get for $x$ large enough:
$$\tilde Z(\Q,C_3\wr C_2;x) \geq Z(K, C_3; x/27) \sim c(K,C_3) x^{1/2}\log(x).$$
For the latter we used the fact that $b(K,C_3)=2$ and that the conjecture is true for abelian groups
\cite{wr}. This already gives a contradiction to Conjecture \ref{con}. Since the asymptotics of the counting
functions of $C_6$ and $S_3(6)$ is $O(x^{1/3})$ we get the contradiction for our group $G$.

Now we introduce a counting function avoiding $\Q(\zeta_3)$.
$$\hat Z(\Q,C_3\wr C_2; x) :=|\left\{L/\Q\mid \exists K\ne
  \Q(\zeta_3):\Gal(L/K)=C_3, [K:\Q]=2, |d_L|\le x\right\}|.$$
Using the averaging results for the $3$-ranks of the class group of
quadratic fields \cite{DaHe} we can prove that
$$\hat Z(\Q,C_3\wr C_2;x) \sim c(C_3 \wr C_2) x^{1/2}$$
as predicted by Malle's conjecture. This means that the cyclotomic intermediate extension
is the reason for the failure of the conjecture.

We remark that all groups of type $C_\ell \wr C_m$, where $\ell>2$
is prime and $\gcd(m,\ell-1)>1$ give counter examples. E.g. let us consider $m=2$.
Then there is a unique quadratic subfield $K=\Q(\sqrt{\pm \ell})$ of $\Q(\zeta_\ell)$.
Since $b(K,C_\ell)=2$ and $b(\Q,C_\ell \wr C_2) =1$ we can derive the same contradiction
as above. We remark that for $\ell>3$ we cannot prove good upper bounds for those groups
since we do not know good estimates for the $\ell$-rank of the class group of quadratic
fields in these cases.

How to fix the conjecture?

One possibility would be to forbid intermediate extensions which are
contained in cyclotomic extensions $\Q(\zeta_\ell)$, where $\ell$ must
be chosen from a set containing all orders of elements of $G$ which
have minimal index. This is more natural in the global function field
setting. Here we restrict to extensions $K/\FF_q(t)$ such that the
normal closure of $K$ contains no constant field extension. Restricting
to those extensions and assuming some heuristic about the
number of points of irreducible varieties over $\FF_q$, Ellenberg and
Venkatesh \cite{ElVe2} are able to deduce the correct upper and lower bounds
in Malle's conjecture including the correct logarithmic power.
We remark that the same type of counter examples apply, if we allow constant
field extensions. E.g. choosing $q\equiv 2 \bmod 3$ and $G=C_3 \wr C_2$ works
as a counter example, when we choose $\FF_{q^2}/\FF_q$ as the quadratic extension.

In order to fix the conjecture we might think that we have to look at maximal
abelian quotients of the given group $G$. Unfortunately this approach is not
sufficient.
E.g. for the group $G=(C_3\wr C_3) \times C_2$ we get $a(G)=1/4, b(\Q,G)=1, b(\Q(\zeta_3),G)=2$.
We can prove that there exist constants $c_1(G), c_2(G)>0$ such that
$$c_1(G)x^{1/4} \leq Z(\Q,G;x) \leq c_2(G)x^{1/4} \;\mbox{ for }x \mbox{ large enough}.$$
Therefore this group does not contradict Conjecture \ref{con}.
Similiar to our original example it is possible to choose $K=\Q(\zeta_3)$ as an intermediate 
extension, but this time it does not change the $\log$-factor.





\end{document}